\documentclass[12pt]{amsart}
\usepackage{graphicx}
\usepackage{amsmath}
\usepackage{amssymb}
\usepackage{amsfonts}
\usepackage{amsthm}
\usepackage{xcolor}

\def\AA{{\mathcal A}}
\def\BB{{\mathcal B}}
\def\FF{{\mathbb F}}
\def\GG{{\mathbb G}}
\def\HH{{\mathbb H}}
\def\LL{{\mathcal L}}
\def\NN{{\mathcal N}}
\def\XX{{\mathcal X}}
\def\YY{{\mathcal Y}}

\def\complex{\mathbb{C}}
\def\CC{\mathcal{C}}
\def\MM{\mathcal{M}}

\def\TT{\mathbb{T}}
\def\real{\mathbb{R}}
\def\integer{\mathbb{Z}}

\def\LOP{{\mathcal N}}
\def\defeq{\displaystyle\mathrel{\mathop=^{\scriptscriptstyle\rm def}}}

\def\id{{\mathbb I}}
\def\eps{{\varepsilon}}
\def\cosi{\mathop{\rm cosi}\nolimits}
\def\even{{\rm e}}
\def\odd{{\rm o}}

\newtheorem{theorem}{Theorem}

\newtheorem{lemma}{Lemma}
\newtheorem{proposition}{Proposition}

\usepackage{hyperref}
\hypersetup{
    colorlinks=true,
    linkcolor=blue,
    filecolor=magenta,
    urlcolor=cyan,
}

\title[Computer assisted proof of branches of solutions]{Computer assisted proof of branches of stationary and periodic solution,
and Hopf bifurcations, for dissipative PDEs}
\date{\today}
\author{Gianni Arioli}
\address{Department of Mathematics, Politecnico di Milano\\
piazza Leonardo da Vinci 32, 20133 Milano, Italy}
\email{gianni.arioli@polimi.it}
\thanks{The author is partially supported by the PRIN project
{\em Equazioni alle derivate parziali di tipo ellittico e parabolico: aspetti geometrici,
disuguaglianze collegate, e applicazioni}.}

\begin{document}
\maketitle
%\tableofcontents
\begin{abstract}
We discuss an approach to the computer assisted proof of the existence of branches
of stationary and periodic solutions for dissipative PDEs,
using the Brussellator system with diffusion and Dirichlet boundary conditions as an example,
We also consider the case where the branch of periodic solutions emanates from 
a branch of stationary solutions through a Hopf bifurcation.
\end{abstract}
\baselineskip11pt

\section{Introduction and main results}
The seminal work of Turing \cite{T} has introduced the explanation of pattern formation in nature through reaction and diffusion of chemical compounds. The Gray-Scott system and the Brussellator with diffusion are well known models of such dynamics.
The literature on these reaction-diffusion systems is very large, see e.g. \cite{AE}-\cite{PPG} and references therein.
Most result on the existence of solutions and bifurcations of branches of solutions for these reaction-diffusion systems concern Neumann boundary conditions, and are based on the explicit knowledge of (constant) stationary solutions.
When other boundary conditions are considered, e.g. homogeneous Dirichlet, no stationary solutions are known explicitly, therefore different techniques must be used to study branches of solutions and bifurcations.

In the last few years different computer assisted methods have been developed to study branches of solutions and bifurcations for both ordinary and partial differential equations, see e.g.  \cite{AK1}-\cite{BQ}. In particular, in \cite{AGK} and \cite{AK2} an efficient method based on the Taylor expansion of the Fourier coefficients of the solution has been introduced. Here we apply such method to study branches of stationary solutions and a Hopf bifurcation for the Bussellator system with diffusion and Dirichlet boundary conditions, and we further expand it to study branches of periodic solutions arising from the Hopf bifurcation, so that the existence of periodic solutions is guaranteed in a given interval of the parameter, instead of some neighborhood of unknown width.

The problem that we choose as an example is the system
\begin{equation}
\begin{cases}
U_t-U_{xx}=\sin x-(b+1)U+U^2V\\
V_t-\frac{1}{64}V_{xx}=bU-U^2V\\
U(0,t)=U(\pi,t)=V(0,t)=V(\pi,t)=0\quad\forall t\,,
\end{cases}
\label{brussel}
\end{equation}
where $b$ is a positive parameter. This is essentially the Brussellator system with diffusion and Dirichlet boundary conditions, the only difference being the $\sin x$ term in the first equation instead of the more conventional constant term used with Neumann boundary conditions (see e.g. \cite{BD}). This choice allows solutions admitting analytic extensions to the whole real line, a property that simplifies the computer assisted estimates, while maintaining the main characteristics of the problem. Note that the diffusion coefficient for $V$ is much smaller than the diffusion coefficient for $U$. Indeed, the bifurcation that we consider here is a typical example of Turing instability, which manifests itself only when the speed of diffusion in the two variables is significantly different.

Our first result is the following:
\begin{theorem}\label{th:stat}
For all $b\in[0,11]$ equation \eqref{brussel} admits an analytic stationary solution $(U_b,V_b)$, symmetric with respect to the reflection $x\mapsto \pi-x$. The map $b\mapsto(U_b,V_b)$ is also analytic. For any fixed value of $b$ the stationary solution is isolated.
\end{theorem}
In particular, Theorem \ref{th:stat} excludes the existence of a symmetry breaking bifurcation in the parameter range considered.

\begin{figure}[hbt!]
\begin{center}
\includegraphics[height=50mm,width=120mm]{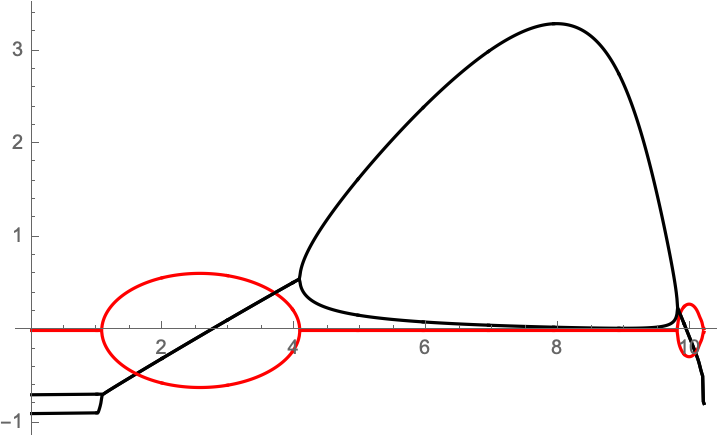}
\caption{\small Eigenvalues.}\label{fig:eigen}
\end{center}
\end{figure}

We studied numerically the linear stability of the solution $(U_{b},V_{b})$ obtained in Theorem \ref{th:stat}. It turns out that at most two eigenvalues of the linearized equation appear to have positive real part. Figure \ref{fig:eigen} represents the real (black) and imaginary (red) parts of such eigenvalues, as $b$ ranges in $[0,11]$. The graph shows that the real part of two complex conjugate eigenvalues changes of sign at $b\simeq2.7$ and $b\simeq 10$, suggesting two supercritical Hopf bifurcations.

Our second results is a proof of the supercritical Hopf bifurcation at $b\simeq2.7$ and of the existence of a branch of periodic solutions emanating from it.

\begin{theorem}\label{th:periodic}
At $b=b_0=2.6992\ldots$ a Hopf bifurcation occurs. More precisely, for all $b\in[b_0,b_1]$ with $b_1=2.7418...$
there exists a periodic solution $(\tilde U_b,\tilde V_b)$, with $\partial \tilde U_b/\partial t\ne0$ for all $b>b_0$, and
$(\tilde U_{b_0},\tilde V_{b_0})=(U_{b_0},V_{b_0})$, where $(U_{b_0},V_{b_0})$ is the solution obtained in Theorem \ref{th:stat}.
Furthermore, the map $b\mapsto(\tilde U_b,\tilde V_b)$ is analytic and the minimal period of the solutions varies in the interval
$[10.35\ldots,10.97\ldots]$.
\end{theorem}

\begin{figure}[hbt!]
\begin{center}
\includegraphics[height=30mm,width=40mm]{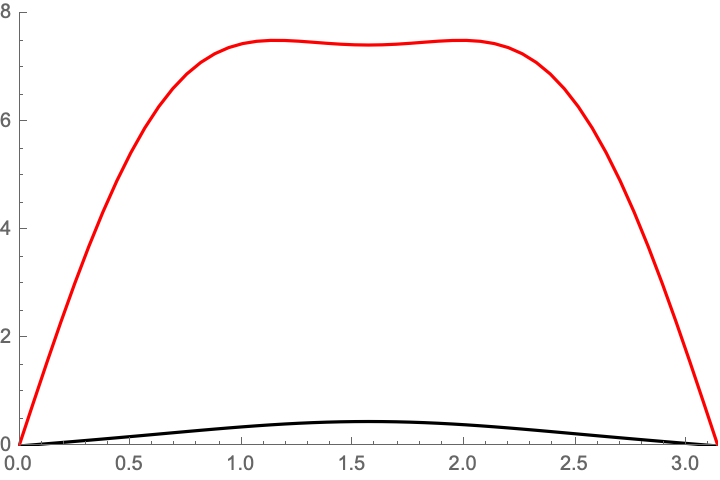}
\includegraphics[height=30mm,width=40mm]{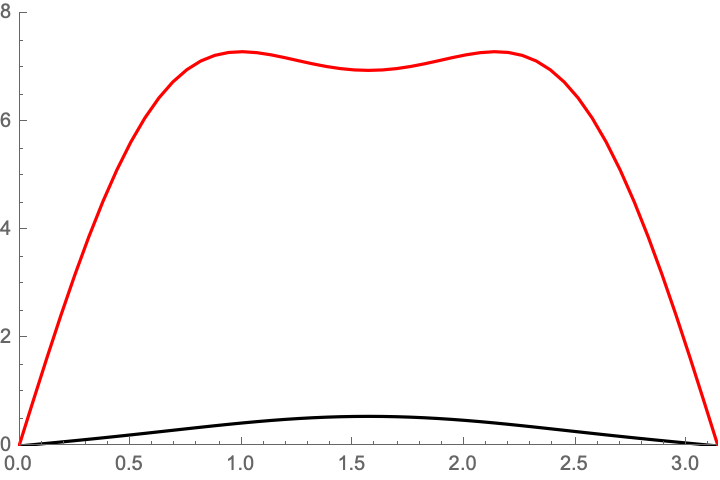}
\includegraphics[height=30mm,width=40mm]{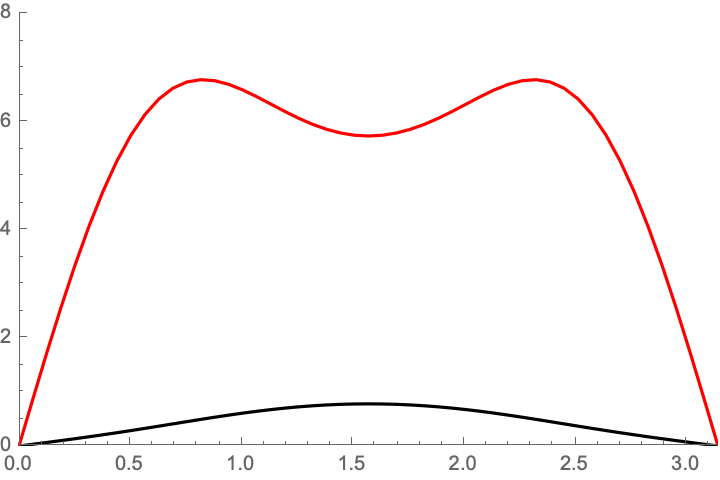}
\includegraphics[height=30mm,width=40mm]{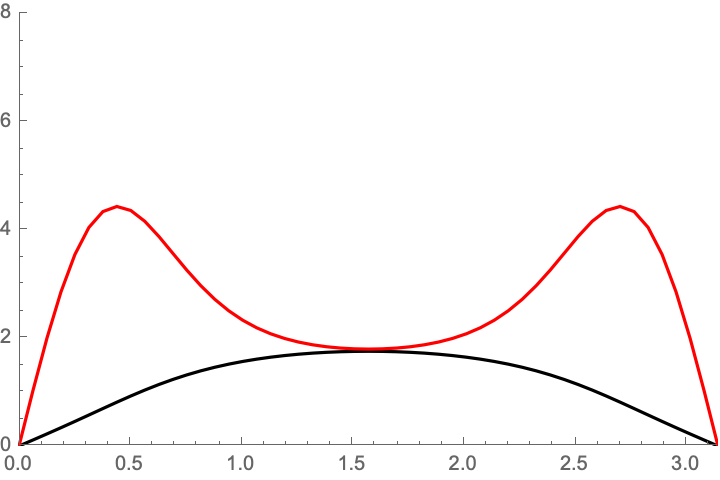}
\includegraphics[height=30mm,width=40mm]{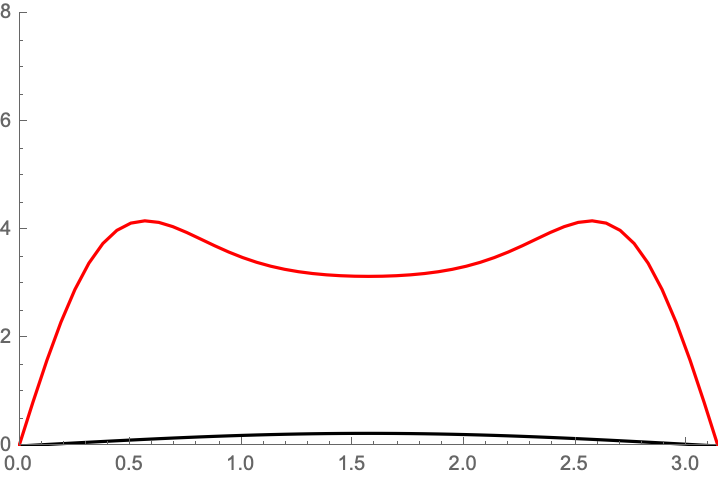}
\includegraphics[height=30mm,width=40mm]{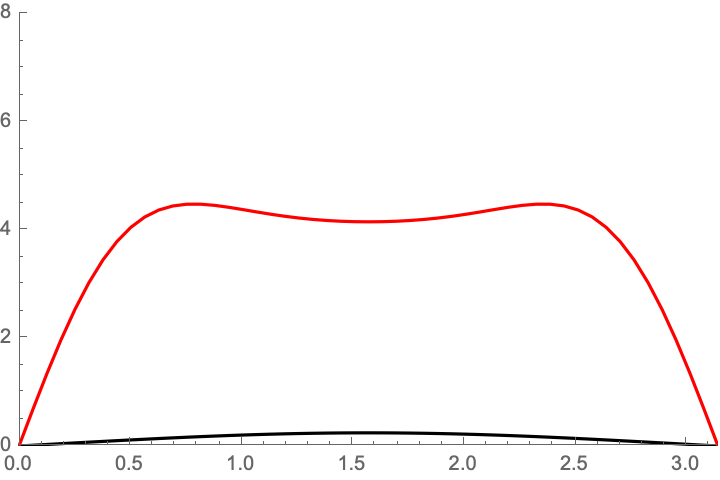}
\includegraphics[height=30mm,width=40mm]{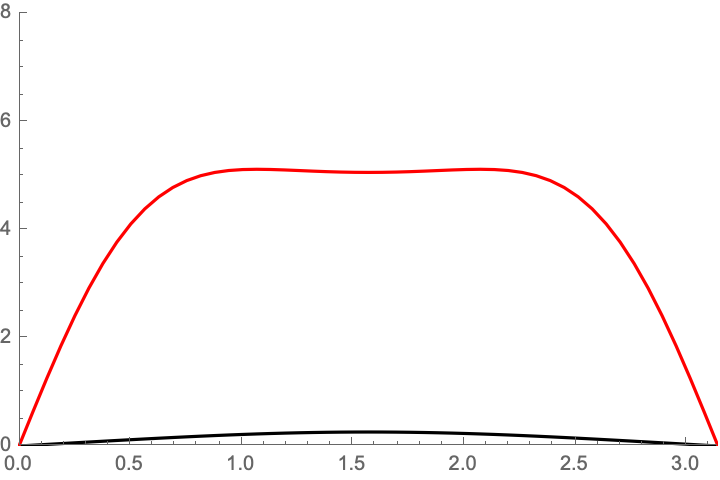}
\includegraphics[height=30mm,width=40mm]{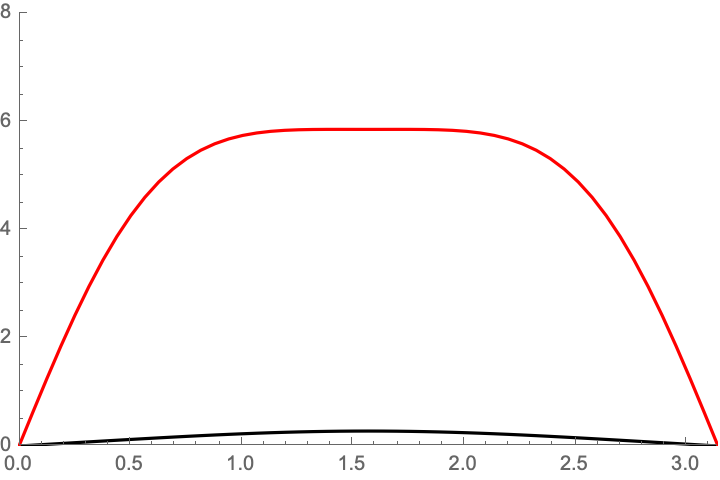}
\includegraphics[height=30mm,width=40mm]{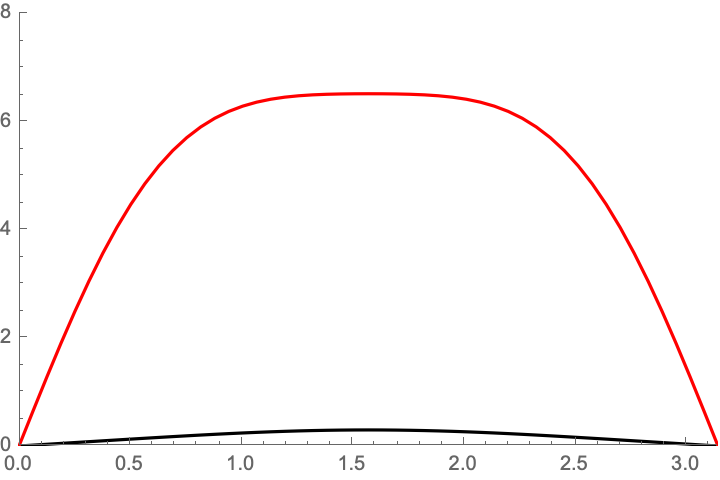}
\includegraphics[height=30mm,width=40mm]{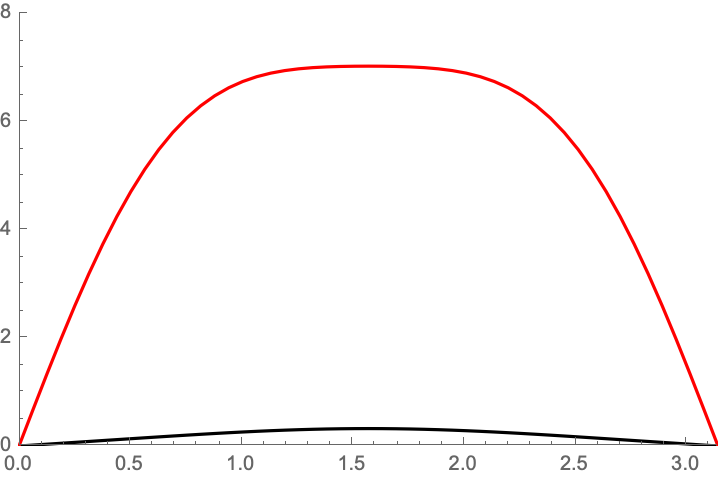}
\includegraphics[height=30mm,width=40mm]{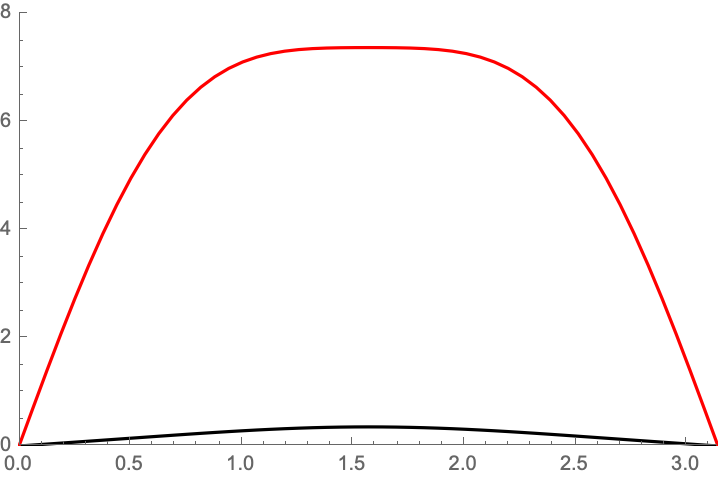}
\includegraphics[height=30mm,width=40mm]{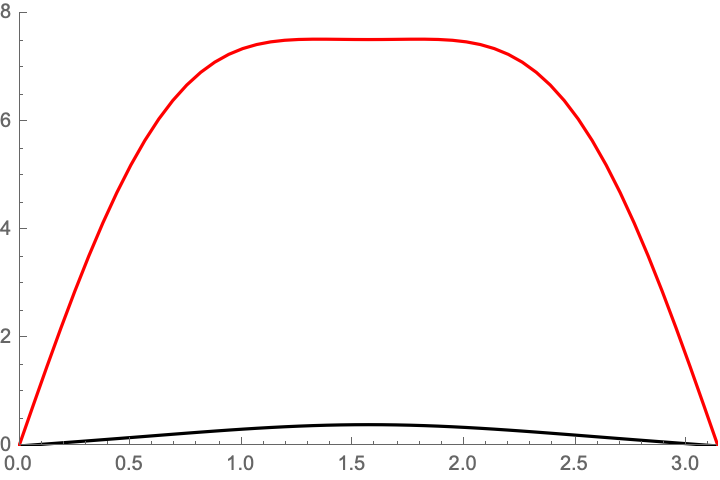}
\caption{\small Snapshots of the periodic solutions at $b=3$ at times $t=kT/12$, $k=0,\ldots,11$. The graph of $u$ is black and the graph of $v$ is red.}\label{fig:periodic}
\end{center}
\end{figure}

\begin{figure}[hbt!]
\begin{center}
\includegraphics[height=40mm,width=50mm]{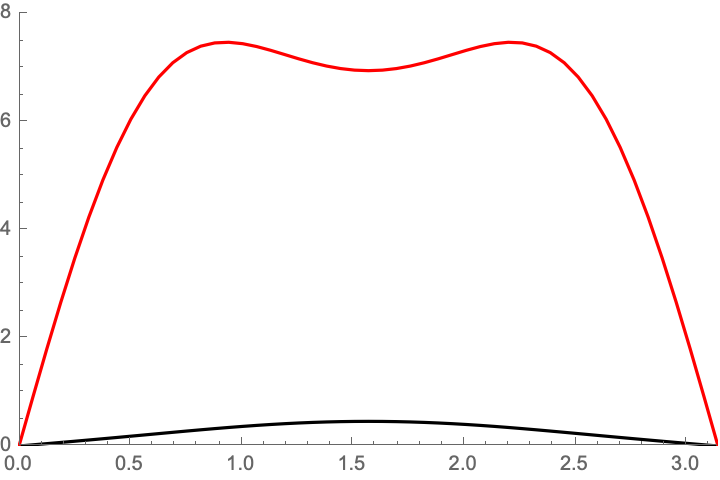}
\hskip1cm
\includegraphics[height=40mm,width=50mm]{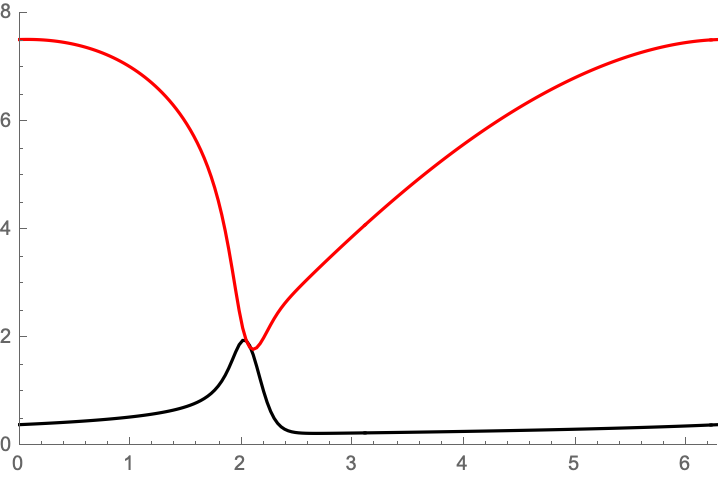}
\caption{\small Stationary solution at $b=3$ (left). Graph of $(\tilde U_3(\pi/2,t),\tilde V_3(\pi/2,t))$ (right).}\label{fig:stat}
\end{center}
\end{figure}

The proof of the existence of the branch of periodic solutions is carried out only in the interval $[b_0,b_1]$ because it would be very computationally expensive to consider a wider interval. As it turns out, the solution varies significantly and the size of the domain of analyticity decreases very fast with increasing $b$, which in turn means that a very large number of Fourier coefficients need to be taken into account. Still, there is no theoretical obstruction to the extension of the branch. For illustration purposes, we also proved the existence of a periodic solution at $b=3$:

\begin{theorem}\label{th:periodic3}
At $b=3$ equation \eqref{brussel} admits an analytic periodic solution $(\tilde U_3,\tilde V_3)$ of period $T=22.91\ldots$.
\end{theorem}
Note that, when $b$ increases from $2.69\ldots$ to $3$, the period of the solution doubles and the graph of the solution changes dramatically. Figure \ref{fig:periodic} shows some snapshots of the solution $(\tilde U_3,\tilde V_3)$, while Figure \ref{fig:stat} shows the stationary solution $(U_3,V_3)$ and the graph of $(\tilde U_3(\pi/2,t),\tilde V_3(\pi/2,t))$.

\section{The strategy of the proofs}
\subsection{Stationary branch}
Let $\rho=65/64$, let $\CC$ be the space of functions
\begin{equation}\label{CCdef}
f=\sum_{k\ge1}f_{k}\sin(kx)\qquad(f_{k}\in\mathbb{C})
\end{equation}
such that
\begin{equation}\label{CCnormdef}
\|f\|:=\sum_{k}|f_{k}|\rho^{k}<+\infty\,.
\end{equation}
Denote by $\CC_S$ the subspace of $\CC$ characterized by the symmetry $f(x)=f(\pi-x)$, that is the subspace of functions with $f_{2k}=0$ for all $k=1,2,\ldots$.
For $u=(U,V)\in\CC^2$ let $\|u\|=\|U\|+\|V\|$. Let $A_b:\CC^2\to\CC^2$ be defined by
\begin{equation}\label{eq:Fstat}
A_b(U,V) := (\partial_{xx})^{-1}(\sin x-(b+1)U+U^2V,64(bU-U^2V))\,,
\end{equation}
Clearly, fixed points of $A_b$ are analytic stationary solutions of \eqref{brussel}. Also note that
$A_b(U,V)\in\CC_S^2$ for all $(U,V)\in\CC_S^2$.

To prove that $A_b$ admits a fixed point $(U_b,V_b)\in\CC_S^2$
for each $b\in[0,11]$, and that the function $b\mapsto(U_b,V_b)$ is analytic,
we write all the coefficients in the Fourier expansion of $(U_b,V_b)$ as Taylor polynomials in $b$:
\begin{align}
(U_b(x),V_b(x))=\sum_{k\text{ odd}} (u_{k}(b),v_k(b))\sin(kx)\,,\\
(u_{k}(b),v_k(b))=\sum_{l=0}^L (u_{kl}(b),v_{kl}(b))\left(b-b_0\over b_1\right)^l\,,
\label{taylor}
\end{align}
where $b_0,b_1$ and $L$ are as in Table \ref{tab:stat}.
As a first step, we choose some Fourier-Taylor
polynomial $\bar u=(\bar U,\bar V)\in\CC_S^2$ that is an approximate fixed point of $A_b$, and some
finite rank operator $M_b:\CC^2\to\CC_S^2$ such that $\id-M_b$ is an approximate inverse of $\id-DA_b(\bar u)$.
Then for $h\in\CC^2$ we define
\begin{equation}
\MM_b(h)=A_b(\bar u+ \Lambda_b h)-\bar u+M_b h\,,
\qquad\Lambda_b=\id-M_b\,.
\label{contraction}
\end{equation}
Clearly, if $h$ is a fixed point of $\MM_b$, then $u=\bar u+\Lambda_b h$ is a fixed point of $A_b$ and, hence,
$u$ solves \eqref{brussel}.
Given $r>0$ and $w\in\CC^2$, let $B_r(w)=\{v\in\CC^2\,:\,\|v-w\|<r\}$.
We partition the interval $[0,11]$ into four subintervals.
The center $b_0$ and width $b_1$ for each subinterval are shown in Table \ref{tab:stat}.
\begin{table}[htp]
\begin{center}
\begin{tabular}{|c|c|c|c|}
\hline
$i$ & $b_0$ & $b_1$ & $L$\\
\hline
1 & 4 & 4 & 36\\
2 & 17/2 & 1/2 & 20\\
3 & 19/2 & 1/2 & 30\\
4 & 21/2 & 1/2 & 30\\
\hline
\end{tabular}
\vskip2mm
\caption{Branch intervals for the stationary solution}
\end{center}
\label{tab:stat}
\end{table}
Then we prove the following lemma with the aid of a computer, see Section \ref{technicalities}.

\begin{lemma}\label{capOne}
The following holds for each $i\in\{1,\ldots,6\}$.
Define $b_0$, $b_1$, and $L$ as in row $i$ of Table \ref{tab:stat}.
There exist a Fourier-Taylor polynomial $\bar u(b)$
of degree $L$  as described in \eqref{taylor},
a bounded linear operator $M_b$ on $\CC$,
and positive real numbers $\eps, r, K$ satisfying $\eps+Kr<r$, such that
\begin{equation}
\|\MM_b(0)\|\le\eps\,,\quad
\|D\MM_b(v)\|\le K\,,\quad \forall v\in B_r(0)
\label{apprfp}
\end{equation}
holds for all $\{b\in\complex: |b-b_0|<b_1\}$. Furthermore, 
$$B_{r^1}(\bar u(b_0^1+b_1^1))\subset B_{r^2}(\bar u(b_0^2-b_1^2))\,,$$
$$B_{r^3}(\bar u(b_0^3+b_1^3))\subset B_{r^4}(\bar u(b_0^4-b_1^4))\,,$$
$$B_{r^3}(\bar u(b_0^3-b_1^3))\subset B_{r^2}(\bar u(b_0^2+b_1^2))\,,$$
where the superscript refers row in Table 1.
\end{lemma}

This lemma, together with the Contraction Mapping Theorem and the Implicit Function Theorem, implies Theorem  \ref{th:stat}.

\subsection{Periodic branch}
The time period $T$ of a (non stationary) periodic solution varies with $b$.
Instead of looking for $T$-periodic solutions of equation \eqref{brussel},
we look for $2\pi$-periodic solutions of
\begin{equation}
\begin{cases}
\alpha U_t-U_{xx}=\sin x-(b+1)U+U^2V\\
\alpha V_t-\frac{1}{64}V_{xx}=bU-U^2V\\
U(0,t)=U(\pi,t)=V(0,t)=V(\pi,t)=0\quad\forall t\,,
\end{cases}
\label{brusselalpha}
\end{equation}
where $\alpha=2\pi/T$ has to be determined.
To simplify notation, a $2\pi$-periodic function
will be identified with a function on the circle $\real/(2\pi\integer)$.

Let $\rho=33/32$ and $\varrho=1+2^{-20}$, let
$$
\cosi_j(t)=\begin{cases}
\cos(jt) &\text{ if }j\ge 0\\
\sin(-jt) &\text{ if }j<0\,,
\end{cases}
$$
and let $\AA$ be the space of functions
$$
f(t)=\sum_{j\in\integer}f_j\cosi_j(t)\,,
$$
with
$$
\|f\|=\sum_{j\in\integer}|f_j|\varrho^{|j|}<+\infty\,.
$$
Given $f\in\AA$, let
$$
f_\even\defeq\sum_{j\text{ even}}f_{j}\cosi_{j}(t)\,,\qquad
f_\odd\defeq\sum_{j\text{ odd}}f_{j}\cosi_{j}(t)\,.
$$
Let $\BB$ be the space of functions
$$
(U(x,t),V(x,t))=\sum_{k\ge1}(U_{k}(t),V_{k}(t))\sin(kx)\,,\\
$$
with $U_k,V_k\in\AA$ and
$$
\|(U,V)\|=\sum_{k\ge1}(\|U_{k}\|+\|V_{k}\|)\rho^k<+\infty\,.
$$

Given $(U,V)\in\BB$, let
$$
(U_\odd,V_\odd)\defeq\sum_{k\ge1}(U_{k\odd}(t),V_{k\odd}(t))\sin(kx)\,,
$$
and define $(U_\even,V_\even)$ similarly.
For $b$ near the bifurcation point $b_0$,
we expect $U,V$ to be nearly time-independent.
So in particular, $(U_\odd,V_\odd)$ is close to zero.
This justifies the following scaling: for some $\beta>0$ define
$$
(U,V)=T_\beta (u,v)\defeq (u_\even+\beta u_\odd,v_\even+\beta v_\odd)\,.
$$
Substituting into \eqref{brussel} yields the equation
\begin{equation}
\begin{cases}
\alpha u_t-d_1u_{xx}=\sin x-(b+1)u+\LOP_s(u,v)\\
\alpha v_t-d_2v_{xx}=bu-\LOP_s(u,v)\,.
\end{cases}
\label{brusselbeta}
\end{equation}
where $s=\beta^2$ and
\begin{equation}\label{eq:lop}
\LOP_s(u,v)
=u_e^2v_e+2u_eu_ov_e+u_e^2v_o+s(u_o^2v_e+2u_eu_ov_o+u_o^2v_o)\,.
\end{equation}
Let
$$\LL_{\alpha,d,b}u\defeq(\alpha \partial_t-d\partial_{xx}+b)^{-1}u\,.$$
Then, if
$$
u(x,t)=\sum_{j,k}u_{jk}\cosi_j(t)\sin{kx}
$$
we have
\begin{equation}\label{eq:LL}
\LL_{\alpha,d,b} u(x,t)=\sum_{k\ge1,j\in\integer}\frac{(dk^2+b)u_{jk}-\alpha ju_{(-j)k}}{(dk^2+b)^2+\alpha^2j^2}\sin(kx)\cosi_j(t)\,,
\end{equation}
and \eqref{brussel} is equivalent to
\begin{equation}
\begin{cases}
u=\FF_s(u,v) := \LL_{\alpha,d_1,b+1}(\sin x+\LOP_s(u,v))\\
v=\GG_s(u,v) := \LL_{\alpha,d_2,0}(bu-\LOP_s(u,v))\,.
\end{cases}
\label{brusselFP}
\end{equation}

One of the features of the equation 
is that  the time-translate of a solution is again a solution.
We eliminate this symmetry by imposing the condition $u_{11}=0$.
In addition, close to the Hopf bifurcation point, we normalize
the odd part of $u$ by choosing $u_{(-1)1}=1$.
This leads to the conditions
\begin{equation}\label{eq:norm}
Au\defeq u_{11}=0\,,\qquad
Bu\defeq u_{(-1)1}=1\,.
\end{equation}
It is convenient to regard $s$ to be the independent parameter
and express $\alpha$ and $b$ as a function of $s$.
The functions $\alpha=\alpha(s)$ and $b=b(s)$
are determined by the condition that $u$
satisfies \eqref{eq:norm}.
Applying $A$ and $B$ to both sides of
$u=\FF_s(u,v)$, and using the identities
$A\partial_{xx}=-A$, $A\partial_t=B$, $B\partial_{xx}=-B$, $B\partial_t=-A$,
we find that
\begin{equation}\label{eq:AB}
\alpha=A\LOP_s(u,v)\,,\qquad b=B\LOP_s(u,v)-1-d_1\,.
\end{equation}
To compute the derivatives of $\{\FF_s,\GG_s\}$, assume that $u,v$ depend on a parameter,
and denote by a dot the derivative with respect to this parameter.
Define
$$
\LL_{\alpha,d,b}'=\partial_t(\alpha\partial_t-d\Delta+b)^{-1}\,.
$$
Then the parameter-derivatives of $\FF_s(u,v)$ and $\GG_s(u,v),$ are given by
\begin{align*}
D\FF_s(u,v)[\dot u,\dot v]&=\LL_{\alpha,d,b+1}(D\LOP_s(u,v)[\dot u,\dot v]-\dot\alpha \LL_{\alpha,d,b+1}'(\sin x+\LOP_s(u,v))\\
&-\dot b\LL_{\alpha,d,b+1}(\sin x+\LOP_s(u,v)))\,,\\
D\GG_s(u,v)[\dot u,\dot v]&=\LL_{\alpha,d,0}(\dot bu+b\dot u-D\LOP_s(u,v)[\dot u,\dot v]\\
&-\dot\alpha \LL_{\alpha,d,0}'(bu-\LOP_s(u,v)))\,,
\end{align*}
where
$$\dot\alpha=AD\LOP_s(u,v)[\dot u,\dot v]\text{ and }\dot b=BD\LOP_s(u,v)[\dot u,\dot v]\,.$$
Away from the Hopf bifurcation point there is no need to normalize
the odd part of $u$, so that $b$ can be a free parameter
and then we choose $\beta=1$.
In this case, imposing the condition $Aw=0$ we find that
$$
\alpha=\frac{A\LOP_1(u,v)}{Bu}\,,
$$
and the parameter-derivatives of $\FF_1(u,v)$ and $\GG_1(u,v),$ are given by
\begin{align*}
D\FF_1(u,v)[\dot u,\dot v]&=\LL_{\alpha,d,b+1}(D\LOP_1(u,v)[\dot u,\dot v]
-\dot\alpha \LL_{\alpha,d,b+1}'(\sin x+\LOP_1(u,v)))\,,\\
D\GG_1(u,v)[\dot u,\dot v]&=\LL_{\alpha,d,0}(b\dot u-D\LOP_1(u,v)[\dot u,\dot v]-\dot\alpha \LL_{\alpha,d,0}'(bu-\LOP_1(u,v)))\,.
\end{align*}

To prove that \eqref{brussel} admits a periodic solution $(U_b,V_b)\in\CC_S^2$
for each $b\in[b_0,b_1]$, and that the function $b\mapsto(U_b,V_b)$ is analytic,
we follow a procedure similar to the stationary case:
we write all the coefficients in the Fourier expansion of $(u,v)$ as Taylor polynomials in $s$:
\begin{equation}
(u_{jk}(s),v_{jk}(s))=\sum_{l=0}^L (u_{jkl},v_{jkl})\left(s-s_0\over \delta\right)^l\,.
\end{equation}
Let
$$
\HH_s(u,v)\defeq(\FF_s(u,v),\GG_s(u,v))\,,
$$
so that periodic solutions of \eqref{brussel} correspond to fixed points of $\HH_s$.
Let $\BB_S\defeq\{(u,v)\in\BB\,:\,u(x,t)=u(\pi-x,t)\,,v(x,t)=v(\pi-x,t)\}$, and note that
$\HH_s(u,v)\in\BB_S$ for all $(u,v)\in\BB_S$.

As a first step, we choose some Fourier-Taylor
polynomial $\bar w=(\bar u,\bar v)\in\BB_S$, that is an approximate fixed point of $\HH_s$, and some
finite rank operator $M=M(s)$ such that $\id-M$ is an approximate inverse of $\id-D\HH_s(\bar w)$
and $Mw\in\BB_S$ for all $w\in\BB_S$.
Then for $h\in\BB_S$ we define
\begin{equation}
\NN_s(h)=\HH_s(\bar w+ \Lambda h)-\bar w+M h\,,
\qquad\Lambda=\id-M\,.
\label{contraction}
\end{equation}
Clearly, if $h$ is a fixed point of $\NN_s$, then $u=\bar u+\Lambda h\in\BB_S$ is a fixed point of $\HH_s$ and, hence,
$u$ solves \eqref{brussel}.
Given $r>0$ and $w\in\BB_S$, let $B_r(w)=\{v\in\BB_S\,:\,\|v-w\|<r\}$.
Lemmas \ref{capTwo} and \ref{capThree} are proved with the aid of a computer, see Section \ref{technicalities}.

\begin{lemma}\label{capTwo}
The following holds for each $i\in\{0,\ldots,5\}$.
Let $s_{0i}=i2^{-10}$, $\delta=2^{-11}$.
There exist a Fourier-Taylor polynomial $\bar w_i(s)$
of degree $5$  as described in \eqref{taylor},
a bounded linear operator $M_i(s)$ on $\CC_S$,
and positive real numbers $\eps_i, r_i, K_i$ satisfying $\eps_i+K_ir_i<r_i$, such that
\begin{equation}
\|\NN_\beta(0)\|\le\eps_i\,,\quad
\|D\NN_\beta(w)\|\le K_i\,,\quad \forall w\in B_{r_i}(0)
\label{apprfp}
\end{equation}
holds for all $\{s\in\complex: |s-s_{0i}|<\delta\}$.
\end{lemma}
Let 
$I$ be the natural embedding of $\CC^2_S$ into $\BB_S$.
\begin{lemma}\label{capThree}
$b(0)=2.6992\ldots$ and $b(11\cdot2^{11})= 2.7418\ldots$ (see \eqref{eq:AB}).
For each $i\in\{0,\ldots,4\}$
$$B_{r_i}(\bar w(s_{0i}+\delta))\subset B_{r^{i+1}}(w(s_{0(i+1)}-\delta))\,,$$
$$B_{r_0}(\bar w(s_0))\subset B_r(I(\bar u(b(0))))\,.$$
\end{lemma}

These lemmas, together with the Contraction Mapping Theorem, imply the following proposition, which in turn yields Theorem \ref{th:periodic}.

\begin{proposition}\label{branchexists}
For each $s\in[0,11\cdot2^{-11}]$ there exists a fixed point $u_s$ of $\HH_s$, $\partial_t u_s\ne0$ for all $s>0$, and the curve $s\mapsto u_s$ is analytic. Furthermore, $u_0$ coincides with the solution $\bar u(b(0))$ of Lemma \ref{capOne}.
\end{proposition}

Finally, Theorem \ref{th:periodic3} follows from the following Lemma, also proved with computer assistance:
\begin{lemma}\label{capFour}
Let $b=3$. There exist a polynomial $\bar w\in\BB_S$
a bounded linear operator $M$ on $\BB_S$,
and positive real numbers $\eps, r, K$ satisfying $\eps+Kr<r$, such that
\begin{equation}
\|\NN_1(0)\|\le\eps\,,\quad
\|D\NN_1(w)\|\le K\,,\quad \forall w\in B_{r}(0)
\label{apprfp}
\end{equation}
holds for all $\{s\in\complex: |s-s_{0i}|<\delta\}$.
\end{lemma}

\section{Technicalities}\label{technicalities}

\subsection{Estimates on the linear operators $\boldmath\LL_{\alpha,d,b}$ and $\boldmath\LL_{\alpha,d,b}'$}
%%%%%%%%%%%%%%%%%%%%%%%%%%%%%%%%%%%%%%%%%%%%%%%%%%%%%%%%%%%%%%%%%%%%%%%%%%%%%%%

Consider the linear operators $\LL_{\alpha,d,b}$
and $\LL_{\alpha,d,b}'$, with $\alpha,d,b\in\real$.
Let $\tilde u=\LL_{\alpha,d,b}u$.
Using \eqref{eq:LL} and the Cauchy-Schwarz inequality in $\real^2$
we have the estimate
\begin{equation}\label{eq:Cjk}
|\tilde u_{j,k}|\le C_{j,k}\sqrt{|u_{j,k}|^2+|u_{j,-k}|^2}\,,\qquad
C_{j,k}=\frac{1}{\sqrt{(dj^2+b)^2+\alpha^2k^2}}\,.
\end{equation}
Note that $C_{j,k}$ is a decreasing function of $j,k$
and can be used to estimate $\LL_{\alpha,d,b}u$
when $u$ is the tail of a Fourier series.

For the operator $\LL_{\alpha,d,b}'$ we have
$$
\tilde u_{jk}
=k\frac{(dj^2+b)u_{jk}-\alpha ku_{j(-k)}}{(dj^2+b)^2+\alpha^2k^2}\,,\qquad
\tilde u=\LL'_{\alpha,d,b}u\,.
$$
A bound analogous to \eqref{eq:Cjk} holds, with
with
$$
C'_{j,k}=\frac{|k|}{\sqrt{(dj^2+b)^2+\alpha^2k^2}}\le \frac{1}{\alpha}\,.
$$

\subsection{The computer assisted part of the proof}

    The methods used here can be considered perturbation theory:
    given an approximate solution, prove bounds
    that guarantee the existence of a true solution nearby.
    But the approximate solutions needed here are too complex
    to be described without the aid of a computer,
    and the number of estimates involved is far too large.

    The first part (finding approximate solutions)
    is a strictly numerical computation.
The rigorous part is still numerical,
but instead of truncating series and ignoring rounding errors,
it produces guaranteed enclosures at every step along the computation.
This part of the proof is written in the
programming language Ada \cite{Ada}.
The following is meant to be a rough guide for the reader who wishes to check the
correctness of our programs.
The complete details can be found in \cite{files}.

In the present context, a ``bound'' on a map $f:\XX\to\YY$
is a function $F$ that assigns to a set $X\subset\XX$
of a given type ({\tt Xtype}) a set $Y\subset\YY$
of a given type ({\tt Ytype}), in such a way that
$y=f(x)$ belongs to $Y$ for all $x\in X$.
In Ada, such a bound $F$ can be implemented by defining
a {\tt procedure F(X\,:\, in Xtype\,;\, Y\,:\, out Ytype)}.

To represent balls in a real Banach algebra $\XX$ with unit ${\bf 1}$
we use pairs {\tt S=(S.C,S.R)},
where {\tt S.C} is a representable number ({\tt Rep})
and {\tt S.R} a nonnegative representable number ({\tt Radius}).
The corresponding ball in $\XX$
is $\langle{\tt S},\XX\rangle=\{x\in\XX:\|x-({\tt S.C}){\bf 1}\|\le{\tt S.R}\}$.

When $\XX=\mathbb{R}$ the data type described above is called {\tt Ball}.
Our bounds on some standard functions involving the type {\tt Ball}
are defined in the packages {\tt Flts\_Std\_Balls}.
Other basic functions are covered in the packages {\tt Vectors} and {\tt Matrices}.
Bounds of this type have been used in many computer-assisted proofs;
so we focus here on the more problem-specific aspects of our programs.

The computation and validation of branches involves Taylor series in one variable,
which are represented by the type {\tt Taylor1} with coefficients of type {\tt Ball}.
The definition of the type and its basic procedures are in the package {\tt Taylors1}.
Given a {\tt Radius} $\rho$, consider the space $\TT_\rho$ of all real analytic functions
$g(t)=\sum_n g_nt^n$ on the interval $|t|<\rho$,
obtained by completing the space of polynomials
with respect to the norm $\|g\|_\rho=\sum_n|g_n|\rho^n$.
Given a positive integer {\tt D},
a {\tt Taylor1} is a triple {\tt P=(P.C,P.F,P.R)},
where {\tt P.F} is a nonnegative integer, ${\tt P.R}=\rho$,
and {\tt P.C} is an {\tt array(0..D) of Ball}.
The corresponding set in $\langle{\tt Taylor1},\TT_\rho\rangle$ is defined as
\begin{equation}
\langle{\tt P},\TT_\rho\rangle
=\sum_{n=0}^{m-1}\bigl\langle{\tt P.C(n)},\real\bigr\rangle p_n
+\sum_{n=m}^D\bigl\langle{\tt P.C(n)},\TT_\rho\bigr\rangle p_n\,,\qquad p_n(t)=t^n\,,
\label{TaylorEnclosure}
\end{equation}
where $m=\min({\tt P.F},{\tt D}+1)$.
For the operations that we need in our proof,
this type of enclosure allows for simple and efficient bounds.

Consider now the space $\AA$
for a fixed domain radius $\varrho>1$ of type {\tt Radius}.
Functions in $\AA$  are represented by the type {\tt Fourier1}
defined in the package {\tt HFouriers1}, which accepts
coefficients in some Banach algebra with unit $\XX$.
In our application the coefficients of {\tt Fourier1} are {\tt Taylor1}.
The type {\tt Fourier1} consists of a triple {\tt F=(F.C,F.E,F.Freq0)},
where {\tt F.C} is an {\tt array(-K,0..K) of Ball},
{\tt F.E} is an {\tt array(-2*K..2*K) of Radius} and {\tt F.Freq0}
is a boolean flag that, when {\tt True}, indicates that the object
{\tt Fourier1} represents a constant function.
The corresponding set $\langle{\tt F},\AA\rangle$
is the set of all function $u=p+h\in\AA$, where
\begin{align*}
p(t)&=\sum_{j=-K}^K \langle{\tt F.C(J)},\XX\rangle\,\cosi_j(t)\,,\qquad h=\sum_{j=-2K}^{2K}h^{j}\,,\\
h^{j}(t)&=\sum_{m\ge0}h^{j}_{2m}\,\cosi_{j+2m}(t)\text{ if }j\ge0\,,\\
h^{j}(t)&=\sum_{m\ge0} h^{j}_{2m}\,\cosi_{j-2m}(t)\text{ if }j<0
\end{align*}
with $\|h^{J}\|\le{\tt F.E(J)}$, for all $J$.
Note that high order errors for the even frequencies and odd frequencies
are handled separately.
For the operations that we need in our proof,
this type of enclosure allows for simple and efficient bounds.
In particular, note that the nonlinearity in \eqref{eq:lop}
requires a nonstandard product, depending on the parameter $s$.
The {\tt procedure Prod} in the package {\tt HFouriers1}
provides bounds for that product.

    We note that {\tt Fourier1} (just like {\tt Taylor1})
    allows a generic type {\tt Scalar} for its coefficients;
    and this {\tt Scalar} can be again a Taylor (or Fourier) series.
    This feature makes it easy to represent Fourier series
    whose coefficients depend on parameters.

To represent the space $\BB$ we use a package {\tt Fouriers1}
which is very similar to the package {\tt HFouriers1}, so that we
do not provide a detailed description here. The main differences consist
in having a standard product, and high order errors for the even frequencies and odd frequencies
are handles together.

More precisely, the package {\tt Fouriers1} is instantiated with scalars of type {\tt Fouriers1}
defined in {\tt HFouriers1}, which in turn has scalars defined in {\tt Taylors1}.
For simplicity, the space $\CC$ is represented with the same object,
with $K=0$.

    A bound on the map $\HH_s$ is implemented
    by the procedure {\tt GMap} in the package {\tt Taylors1.Foor.Fix}.
    Defining and estimating a contraction like $\NN_\beta$
    is a common task in many of our computer-assisted proofs.
    An implementation is done via two generic packages,
    {\tt Linear} and {\tt Linear.Contr}.
    For a description of this process we refer to \cite{AKiii}.

\end{document}